\documentclass[showkeys,12pt,times,article,onecolumn,superscriptaddress]{revtex4-2}
\usepackage{amssymb}
\usepackage{amsmath}
\usepackage{amsthm}
\usepackage[utf8]{inputenc}
\usepackage[english]{babel}
\usepackage{color}
\usepackage{comment}
\usepackage{geometry}		
\usepackage{graphicx}
\usepackage{epstopdf}
\usepackage{booktabs}
\usepackage[colorlinks=true,
            linkcolor=red,
            urlcolor=red,
            citecolor=red]{hyperref}
\usepackage[mathlines]{lineno} 
\usepackage{amsmath, amsfonts}
\usepackage{physics}
\usepackage{booktabs}
\usepackage{tabularx}
\usepackage{subfigure}
\usepackage{float}
\usepackage{amscd}

\usepackage{amsfonts}
\usepackage{bbding}
\usepackage{float}
\usepackage{mathrsfs}
\usepackage[T1]{fontenc}
\usepackage{lineno}
\usepackage{xcolor}

\let\originalleft\left
\let\originalright\right
\renewcommand{\left}{\mathopen{}\mathclose\bgroup\originalleft}
\renewcommand{\right}{\aftergroup\egroup\originalright}

\newtheorem{theorem}{Theorem}

\newtheorem{corollary}[theorem]{Corollary}
\newtheorem{lemma}[theorem]{Lemma}
\newtheorem{proposition}[theorem]{Proposition}
\theoremstyle{definition}

\theoremstyle{remark}

\begin{document}
\title{
Limit Cycles Appearing from a Generalized Heteroclinic Loop \\
with an Elementary Saddle and a Nilpotent Saddle \\
by Perturbing Piecewise Hamiltonian Systems
}

\author{Zhou Jin}
\affiliation{School of Mathematics and Physics, China University of Geoscience, Wuhan 430074, P. R. China}

\author{Zhouchao Wei}
\affiliation{School of Mathematics and Physics, China University of Geosciences, Wuhan 430074, P. R. China}

\author{Sishu Shankar Muni}
\affiliation{
Department of Physical Sciences, Indian Institute of Science Education and Research Kolkata, India}
\maketitle
\section*{abstract}
In this paper, we study limit  cycle bifurcations for a class of general near-Hamiltonian systems near a heteroclinic loop with an elementary saddle and a nilpotent saddle. Firstly, we consider the behaviors of the unperturbed system, providing the phase portraits of the system and the necessary conditions for the appearance of a  heteroclinic loop with an elementary saddle and a nilpotent saddle by using the relevant qualitative theory. Then, with consideration of the expression of the first-order Melnikov function, we derive its expansion near the heteroclinic loop by employing some techniques and properties of Abelian integral. Finally, we investigate the coefficients of the expansion, and show that there can exist at least $4[\frac{n+1}{2}]+1$ limit cycles under disturbance.

\section{ Introduction and Main Results}
\noindent It is well-known that the second part of Hilbert's 16th problem, which is to investigate the maximal number of limit cycles that a planar polynomial differential system can have, still remains open (see \cite{hilbert1902mathematical}). In 1977, Arnold proposed the so-called weakened Hilbert's 16th problem, which is to study the maximal number of simple zeros of the Abelian integral (see \cite{arnold1977loss}). Recently, extensive studies on piecewise dynamical systems have caught attention of researchers (see \cite{castillo2017pseudo, llibre2019limit, wang2019number, chen2022melnikov}). This class of systems plays an important role and has lots of applications in many fields including biology \cite{Sohail2022, Zeb2022, Xu2022}, engineering \cite{Miao2020, Zhao2022}, neural network \cite{Aouiti2019, Liu2020, He2022}, control theory \cite{Jedda2019, Hu2020, Zhang2022}.

	Piecewise  dynamical systems exhibit more complex dynamical behavior than the smooth cases, and there are still a large number of unsolved problems that is yet unexplored. To the best of our knowledge, enormous results focus on the maximal number of limit cycles appearing in piecewise differential systems under the condition that the plane is divided into two parts by one switching manifold with small perturbations so far, see, e.g. \cite{Wei2015, Sui2019, Chen2020, Wang2021, Ji2022, Ghermoul2022}. Until now, there exist some distinguished results, which are concerned on estimates of the number of limit cycles near a heteroclinic orbit. In 2018, the number of limit cycles was derived  for piecewise differential systems having a generalized heteroclinic loop with an elementary saddle and a nilpotent saddle or with an elementary saddle and a cusp \cite{Liu2018}. In 2021, it was proved that there are at least $3n-1$ limit cycles emerging from a generalized heteroclinic loop with a cusp and a nilpotent saddle in a piecewise cubic polynomial system \cite{Xiong2021}.
Now we consider the following class of piecewise smooth systems:
\begin{equation}\label{eq_g}
\begin{aligned}
&
\begin{cases}
	\dot{x}=y, \\
	\dot{y}=a(x-b)^r,  \\
\end{cases}
x<0,\\
&
\begin{cases}
	\dot{x}=y, \\
	\dot{y}=c(x-d)^s,
\end{cases}
x\geq 0,
\end{aligned}
\end{equation}
where $a> 0$, $b<0$, $c>0$, $d>0$, $r\geq 1$, $s\geq 1$ and $\dfrac{1}{r+1}a(-b)^{r+1}=\dfrac{1}{s+1}c(-d)^{s+1}$. 
In particular, when $s=1$, $r$ odd (even) and $r\geq 2$,  system (\ref{eq_g}) will be the considered system in \cite{Liu2018}. When $r=2$ and $s=3$,  system (\ref{eq_g}) will be the considered system in  \cite{Xiong2021}.

In this paper, by using the first order Melnikov function, we investigate the number of limit cycles bifurcating from a generalized heteroclinic loop with an elementary saddle and a nilpotent saddle of a switching differential system described by
\begin{equation}\label{eq1}
\begin{aligned}
&
\begin{cases}
	\dot{x}=y, \\
	\dot{y}=a(x-b)^3,  \\
\end{cases}
x<0,~ \\
&
\begin{cases}
	\dot{x}=y, \\
	\dot{y}=c(x-d),
\end{cases}
x\geq 0,
\end{aligned}
\end{equation}
where
\begin{equation}\label{condition}
	a>0,~ b<0,~ c>0,~ d>0,~ \frac{1}{4}ab^4=\frac{1}{2}cd^2.
\end{equation}
Then we perturb system~(\ref{eq1}) with piecewise polynomials of degree $n\geq 1$. It follows that
\begin{equation}\label{eq_pw}
\begin{aligned}
&
\begin{cases}
	\dot{x}=y+\varepsilon p^-(x,y), \\
	\dot{y}=a(x-b)^3+\varepsilon q^-(x,y), \\
\end{cases}
x<0, \\
&
\begin{cases}
	\dot{x}=y+\varepsilon p^+(x,y), \\
	\dot{y}=c(x-d)+\varepsilon q^+(x,y),
\end{cases}
x\geq 0,
\end{aligned}
\end{equation}
where $p^\pm(x,y)=\sum\limits_{i+j=0}^n a_{ij}^\pm x^iy^j$ and $q^\pm(x,y)=\sum\limits_{i+j=0}^n b_{ij}^\pm x^iy^j$ are polynomials of degree $n$, $\varepsilon >0$ is a small parameter, $a_{ij}^\pm$ and $b_{ij}^\pm$ are arbitrary coefficients.

Employing Melnikov function method, we have the following theorem.
\begin{theorem}\label{thm1}
Assume that~(\ref{condition}) holds. For system~(\ref{eq_pw}), the maximal number of limit cycles near the generalized heteroclinic loop $L$, which is equal to the maximal number of simple zeros of $I(h)$ for $h$ near $\dfrac{1}{4}ab^4$, is $4[\dfrac{n+1}{2}]+1$ and it can reach the bound.
\end{theorem}

According to Theorem 4 in \cite{Liu2018}, if it takes $m=3$, we can get $n+2[\dfrac{n-1}{2}]+1$ limit cycles for $n<2$, $n+2[\dfrac{n-1}{2}]-[\dfrac{2n-4}{4}]$ limit cycles for $2\leq n\leq 4$ and $n+[2n]-\dfrac{9}{4}+2-\displaystyle\sum\limits_{i=0}^{[\frac{n-2}{4}]}[\dfrac{n}{2}-(2i+1)]$ limit cycles for $n>4$ near $L$ in system (\ref{eq_pw}). As a comparison, we get more limit cycles when $n\leq 4$ or $n\geq 15$ by simple verification. Then, compared to the result in \cite{Xiong2021}, where the maximal number of limit cycles near a heteroclinic loop is $3n-1$ for taking $r=2$ and $s=3$ in system (\ref{eq_g}), we have more limit cycles when $n=1,3$ and same number of limit cycles when $n=2$, which means we can get more limit cycles near a heteroclinic loop in a simpler system by applying some low-order polynomial perturbations. The proof of Theorem \ref{thm1} is given in section \ref{proof}.
\section{Heteroclinic Bifurcations of System (\ref{eq_pw})}
\begin{figure}
	\centering
	\includegraphics[angle=0,width=12.0cm]{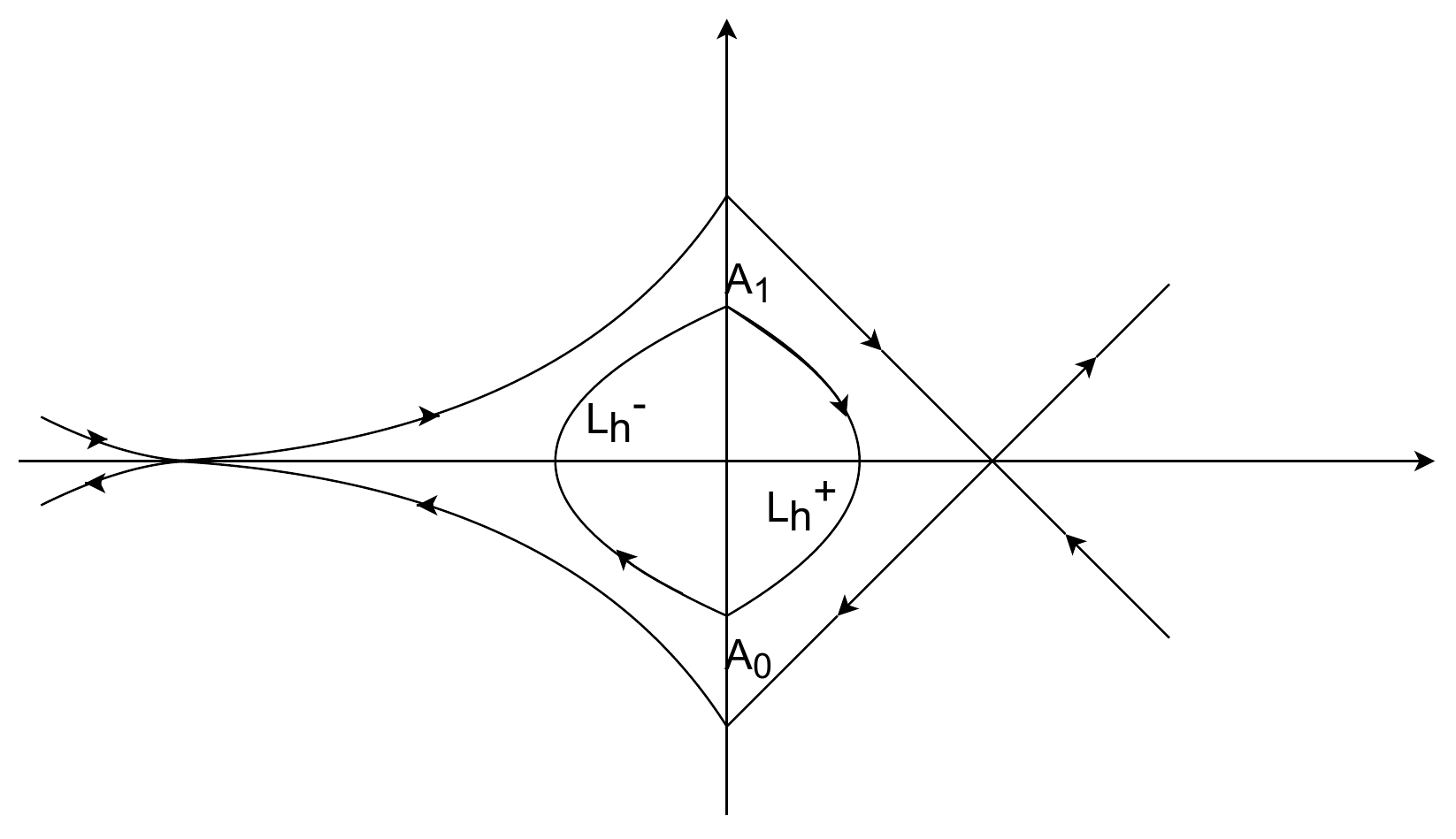}\\
	\caption{Phase portrait of system~(\ref{eq1}).}
	\label{fig1}
\end{figure}

It can be shown that the corresponding Hamiltonian functions  for system~(\ref{eq1}) are given by
\begin{equation}\label{eq_H}
H(x,y)=
\begin{cases}
H^-(x,y)=\frac{1}{2}y^2-\frac{1}{4}ax^4+abx^3-\frac{3}{2}ab^2x^2+ab^3x,~x<0,\\
H^+(x,y)=\frac{1}{2}y^2-\frac{1}{2}cx^2+cdx,~x\geq 0.
\end{cases}
\end{equation}
System~(\ref{eq1}) has a generalized elementary center at $(0,0)$ (the definition about `elementary center' can be found in \cite{Han2010}). There is a family of periodic orbits from (\ref{eq_H}) given by
\begin{equation}\label{eq_L_h}
\begin{aligned}
L_h &= L_h^-\cup L_h^+\\
&= \{(x,y)\,|\,H^-(x,y)=h,x<0\}\cup\{(x,y)\,|\,H^+(x,y)=h,x\geq 0\},\, h\in (0,\frac{1}{4}ab^4),
\end{aligned}
\end{equation}
which intersect the $y$-axis at $A_1(h)=(0,\sqrt{2h})$ and $A_0(h)=(0,-\sqrt{2h})$. As $h\rightarrow 0$, $L_h$ approaches to $(0,0)$. As $h\rightarrow\dfrac{1}{4}ab^4$, $L_h$ approaches to a heteroclinic loop $L$ with a saddle point $(d,0)$ and a nilpotent saddle point $(b,0)$ (see Fig.~\ref{fig1}).
By Theorem 1.1 of \cite{LIU2010}, the generalized first order Melnikov function of system~(\ref{eq_pw}) can be written as
\begin{equation}
	I(h)=I^-(h)+I^+(h),\quad h\in (0,\frac{1}{4}ab^4),
\end{equation}
where
\begin{equation}\label{M_h}
I^{\pm}(h)=\int\limits_{L_h^{\pm}}\sum_{i+j=0}^nb_{ij}^{\pm}x^iy^jdx-\sum_{i+j=0}^na_{ij}^{\pm}x^iy^jdy.	
\end{equation}

\subsection{The Algebraic Structure of $I(h)$}

At the beginning of this section, we introduce some integrals for convenience, which are given by
\begin{equation}\label{I_i}
\begin{aligned}
	C_{ik}^-(h)&=\int\limits_{(\frac{a}{2})^{\frac{1}{4}}b}^{-(\frac{ab^4}{2}-2h)^{\frac{1}{4}}} t^i(2h+t^4-\frac{1}{2}ab^4)^{k+\frac{1}{2}}dt,\\
	C_{ik}^+(h)&=\int\limits_{(cd^2-2h)^{\frac{1}{2}}}^{c^{\frac{1}{2}}d}t^i(2h+t^2-cd^2)^{k+\frac{1}{2}}dt,\quad  
\end{aligned}
\end{equation}
where $h\in(0,\frac{1}{4}ab^4), and i=0,1,\cdots$.
Then, we obtain the following lemma.
\begin{lemma}\label{lemma1}
	For $h\in (0,\dfrac{1}{4}ab^4)$, we have the following expansions:
	\begin{equation*}
	\begin{aligned}
	& C_{ik}^-(h)
	= 
	\begin{cases}
	\alpha_i(1)(2h)^{{\frac{3}{2}}}+\beta_i(1)(2h-\frac{1}{2}ab^4)C_{i0}^-(h),\quad k=1,\\
	\alpha_i(k)(2h)^{k+\frac{1}{2}}+\sqrt{2h}\sum\limits_{j=1}^{k-1}\alpha_i(k-j)\prod\limits_{l=0}^{j-1}\beta_i(k-l)(2h-\frac{1}{2}ab^4)^j(2h)^{k-j}\\
	    \quad +\prod\limits_{j=1}^k\beta_i(j)(2h-\frac{1}{2}ab^4)^kC_{i0}^-(h),\quad k\geq 2,
	\end{cases}
	\end{aligned}
	\end{equation*}
and
\begin{equation*}
\begin{aligned}
& C_{ik}^+(h)
= 
\begin{cases}
\gamma_i(1)(2h)^{\frac{3}{2}}+\delta_i(1)(2h-cd^2)C_{i0}^+(h),\quad k=1,	\\
\gamma_i(k)(2h)^{k+\frac{1}{2}}+\sqrt{2h}\sum\limits_{j=1}^{k-1}\gamma_i(k-j)\prod\limits_{l=0}^{j-1}\delta_i(k-l)(2h-cd^2)^j(2h)^{k-j}\\
\quad +\prod\limits_{j=1}^{k}\delta_i(j)(2h-cd^2)^kC_{i0}^+(h),\quad k\geq 2,
\end{cases}	
\end{aligned}
\end{equation*}
where $i=0,1,\cdots$ and
\begin{equation}\label{eq_kappa}
\begin{aligned}
	&\alpha_i(k)=\frac{-1}{i+1+4(k+\frac{1}{2})}(\frac{a}{2})^{\frac{i+1}{4}}b^{i+1},
	&\beta_i(k)=\frac{4(k+\frac{1}{2})}{i+1+4(k+\frac{1}{2})},\\
	&\gamma_i(k)=\frac{1}{i+1+2(k+\frac{1}{2})}c^{\frac{i+1}{2}}d^{i+1},
	&\delta_i(k)=\frac{2(k+\frac{1}{2})}{i+1+2(k+\frac{1}{2})},
\end{aligned}
k\in\mathbb{N}^+.
\end{equation}
\end{lemma}
\begin{proof}
	Firstly, we have
	\begin{equation*}
	\begin{aligned}
	\int t^i(2h+t^4-\frac{1}{2}ab^4)^{k+\frac{1}{2}}dt & =\frac{1}{i+1+4(k+\frac{1}{2})}\int(\frac{2h-\frac{1}{2}ab^4}{t^4}+1)^{k+\frac{1}{2}}d(t^{i+1+4(k+\frac{1}{2})}) \\
	& =\frac{1}{i+1+4(k+\frac{1}{2})}t^{i+1}(2h+t^4-\frac{1}{2}ab^4)^{k+\frac{1}{2}} \\
	& +\frac{4(2h-\frac{1}{2}ab^4)(k+\frac{1}{2})}{i+1+4(k+\frac{1}{2})}\int t^i(2h+t^4-\frac{1}{2}ab^4)^{k+\frac{1}{2}-1}dt, \\
	\int t^i(2h+t^2-cd^2)^{k+\frac{1}{2}}dt & =\frac{1}{i+1+2(k+\frac{1}{2})}\int (\frac{2h-cd^2}{t^2}+1)^{k+\frac{1}{2}}d(t^{i+1+2(k+\frac{1}{2})}) \\
	& =\frac{1}{i+1+2(k+\frac{1}{2})}t^{i+1}(2h+t^2-cd^2)^{k+\frac{1}{2}} \\
	& +\frac{2(2h-cd^2)(k+\frac{1}{2})}{i+1+2(k+\frac{1}{2})}\int t^i(2h+t^2-cd^2)^{k+\frac{1}{2}-1}dt.
\end{aligned}
\end{equation*}
Then, we get
\begin{equation*}
\begin{aligned}
	& C_{ik}^-(h)  =\alpha_i(k)(2h)^{k+\frac{1}{2}}+\beta_i(k)(2h-\frac{1}{2}ab^4)C_{i,k-1}^-(h), \\
	& C_{ik}^+(h)  =\gamma_i(k)(2h)^{k+\frac{1}{2}}+\delta_i(k)(2h-cd^2)C_{i,k-1}^+(h),
\end{aligned}
\end{equation*}
So Lemma \ref{lemma1} holds for all $k\in\mathbb{N}^+$ by induction.
\end{proof}

\begin{lemma}\label{lemma2}
	Assume that~(\ref{eq_H}) and~(\ref{eq_L_h}) hold, then $I^\pm (h)$ will have the following forms:
\begin{equation}\label{eq_M(h)1}
\begin{aligned}
I^+(h) & =\sum\limits_{i+2k=0}^{n-1}A_{i,2k}^+I_{i,2k}^+(h)+\sum\limits_{k=0}^{[\frac{n}{2}]}\frac{2a_{0,2k}^+}{2k+1}(2h)^{k+\frac{1}{2}}, \\
I^-(h) & =\sum\limits_{i+2k=0}^{n-1}A_{i,2k}^-I_{i,2k}^-(h)-\sum\limits_{k=0}^{[\frac{n}{2}]}\frac{2a_{0,2k}^-}{2k+1}(2h)^{k+\frac{1}{2}},
\end{aligned}
h\in (0,\frac{1}{4}ab^4),
\end{equation}
where
\begin{equation}\label{eq_A_i,2k}
A_{i,j}^+=b_{i,j+1}^++\frac{i+1}{j+1}a_{i+1,j}^+,
\end{equation}
\begin{equation}\label{eq_I_i,2k}
	I_{i,j}^\pm(h)=\int\limits_{L_h^\pm}x^iy^{j+1}dx.       
\end{equation}
Here $A_{i,2k}$, $i+2k=0$, $1$ and $\cdots$, $n-1$, $a_{0,2k}^+$, $k=0$, $1$, $\cdots$ $[\frac{n}{2}]$ are independent with each other.
\end{lemma}
\begin{proof}
The first-order Melnikov function $I(h)$ has the following expansion:
\begin{equation}
I^+(h)=\sum_{i+j=0}^nb_{ij}^+\int\limits_{L_h^+}x^iy^jdx-\sum_{i+j=0}^n\oint\limits_{L_h^+\cup\overrightarrow{A_0A_1}}a_{ij}^+x^iy^jdy+\sum_{i+j=0}^n\int\limits_{\overrightarrow{A_0A_1}}a_{ij}^+x^iy^jdy.
\end{equation}
Hence, by using Green formulas, one obtains that
\begin{equation}
I^+(h)=\sum_{i+j=0}^{n-1}(b_{i,j+1}^++\frac{i+1}{j+1}a_{i+1,j}^+)\int\limits_{L_h^+}x^iy^{j+1}dx+\sum_{k=0}^{[\frac{n}{2}]}\frac{2a_{0,2k}^+}{2k+1}(2h)^{k+\frac{1}{2}}.
\end{equation}
When $j$ is odd, let $j=2k+1, k\geq 0$, one has
\begin{equation*}
\int\limits_{L_h^+}x^iy^{j+1}dx=\int\limits_{L_h^+}x^iy^{2k+2}dx=\int\limits_{L_h^+}x^i(2h+cx^2-cdx)^{k+1}dx\equiv 0.
\end{equation*}
Summarizing the above discussions, it gives the form of $I^+(h)$ in~(\ref{eq_M(h)1}).
Similarily, we can get the expression of $I^-(h)$.
\end{proof}

Now, we can denote $I^\pm (h)$ as
\begin{equation}
I^\pm (h)=I_n^\pm (h)+\sum_{k=0}^{[\frac{n}{2}]}\frac{2a_{0,2k}^+}{2k+1}(2h)^{k+\frac{1}{2}},~
h\in(0,\frac{1}{4}ab^4),
\end{equation}
where
\begin{equation}
I_n^\pm (h)=\sum_{i+2k=0}^{n-1}A_{i,2k}^\pm I_{i,2k}^\pm (h).
\end{equation}
Then, we have the following results.
\begin{lemma}\label{lemma3}
For system~(\ref{eq_pw}), $I_n^\pm (h)$ have the following expansions:
\begin{equation}\label{eq_M_n(h)}
\begin{aligned}
	& I_n^+(h)=\sum\limits_{k=0}^{[\frac{n-1}{2}]}B_{0,2k}^+I_{0,2k}^+(h)+\sum\limits_{k=0}^{[\frac{n-2}{2}]}B_{1,2k}^+(2h)^{k+\frac{3}{2}},~n\geq 2, \\
	& I_n^-(h)=\sum\limits_{k=0}^{[\frac{n-1}{2}]}B_{0,2k}^+I_{0,2k}^-(h)+\sum\limits_{k=0}^{[\frac{n-2}{2}]}B_{1,2k}^-I_{1,2k}^-(h)+\sum\limits_{k=0}^{[\frac{n-3}{2}]}B_{1,2k}^-I_{2,2k}^-(h)+\sum\limits_{k=0}^{[\frac{n-4}{2}]}B_{3,2k}^-(2h)^{k+\frac{3}{2}},~n\geq 4,
	\end{aligned}
\end{equation}
where
\begin{equation}\label{barA}
\begin{aligned}
	& B_{0,2k}^+=A_{0,2k}^++L(A_{i,2k}^+|i+2k\leq n-1,~i\geq 1,~k\geq0),\\
	& B_{1,2k}^+=L(A_{i,2k}^+|i+2k\leq n-1,~i\geq 1,~k\geq0),\\
	& B_{1,2k}^-=A_{1,2k}^-+L(A_{i,2k}^-|i+2k\leq n-1,~i\geq 3,~k\geq0),\\
	& B_{2,2k}^-=A_{2,2k}^-+L(A_{i,2k}^-|i+2k\leq n-1,~i\geq 3,~k\geq0),\\
	& B_{3,2k}^-=A_{3,2k}^-+L(A_{i,2k}^-|i+2k\leq n-1,~i\geq 3,~k\geq0),\\
	& B_{4,2k}^-=L(A_{i,2k}^-|i+2k\leq n-1,~i\geq 3,~k\geq0)
\end{aligned}
\end{equation}
and $L(\bullet )$ denotes a linear combination.
\end{lemma}

\begin{proof}
	Firstly, differentialing both sides of $H^+(x,y)=h$, we have
\begin{equation*}
cxdx-cddx-ydy=0.
\end{equation*}
Multiplying both sides by $x^{i-1}y^{2k+1}$ and integrating along $L_h^+$, one has
\begin{equation}\label{eq_integral}
	cI_{i,2k}^+(h)-cdI_{i-1,2k}^+(h)-\int\limits_{L_h^+}x^{i-1}y^{2k+2}dy=0.
\end{equation}
Then, we have
\begin{equation*}
\begin{aligned}
	\int\limits_{L_h^+}x^{i-1}y^{2k+2}dy & =-\oint\limits_{L_h^+\cup\overrightarrow{A_0A_1}}\frac{i-1}{2k+3}x^{i-2}y^{2k+3}dx-\int\limits_{\overrightarrow{A_0A_1}}x^{i-1}y^{2k+2}dy  =
    \begin{cases}
    \begin{aligned}
    	& -\int\limits_{\overrightarrow{A_0A_1}}y^{2k+2}dy, & i=1, \\
    	& -\frac{i-1}{2k+3}I_{i-2,2k+2}^+(h), & i \geq 2.
    \end{aligned}
    \end{cases}
\end{aligned}
\end{equation*}
Then, with~(\ref{I_i}) and ~(\ref{eq_integral}), it implies that
\begin{equation}\label{eq_I_i,2k1}
\begin{aligned}
& I_{i,2k}^+(h)=dI_{i-1,2k}^+(h)-\frac{i-1}{2k+3}\frac{1}{c}I^+_{i-2,2k+2}, \\
& I_{1,2k}^+(h)=dI_{0,2k}^+(h)-\frac{1}{c}\int\limits_{\overrightarrow{A_0A_1}}y^{2k+2}dy=dI_{0,2k}^+(h)-\frac{2}{2k+3}\frac{1}{c}(2h)^{k+\frac{3}{2}},
\end{aligned}
\end{equation}
where $i\geq 2$, $k\geq 0$.

Further, we can prove the form of $I_n^+(h)$ in~(\ref{eq_M_n(h)}) by using mathematical induction. Firstly, let $n=2$, one has that
\begin{equation}
\begin{aligned}
	I_2^+(h)& =A_{0,0}^+I_{0,0}^+(h)+A_{1,0}^+I_{1,0}^+(h) \\
	& =A_{0,0}^+I_{0,0}^+(h)+A_{1,0}^+(dI_{0,0}^+(h)-\frac{2}{3}\frac{1}{c}(2h)^\frac{3}{2}) \\
	& =(A_{2,0}^++dA_{1,0}^+)I_{0,0}^+(h)-\frac{2}{3c}A_{1,0}^+(2h)^\frac{3}{2},
\end{aligned}
\end{equation}
which means~(\ref{eq_M_n(h)}) holds for $n=2$. Next, assuming the conclusion holds for $n$, $n\geq 2$, it is easy to find that
\begin{equation*}
\begin{aligned}
	I_{n+1}^+(h) & =\sum\limits_{i+2k=0}^nA_{i,2k}^+I_{i,2k}^+(h) \\
	& =\sum\limits_{k=0}^{[\frac{n}{2}]}A_{0,2k}^+I_{0,2k}^+(h)+\sum\limits_{k=0}^{[\frac{n-1}{2}]}A_{1,2k}^+I_{i,2k}^+(h)+\sum\limits_{i+2k=0,i\geq 2}^nA_{1,2k}^+I_{i,2k}^+(h).
\end{aligned}
\end{equation*}
From~(\ref{eq_I_i,2k1}), it implies that
\begin{equation*}
\begin{aligned}
	& \sum\limits_{i+2k=0,i\geq 2}^{n}A_{i,2k}^+I_{i,2k}^+(h)\\
	= & \sum\limits_{i+2k=0,i\geq 2}^n A_{i,2k}^+(dI_{i-1,2k}^+(h)-\frac{i-1}{2k+3}\frac{1}{c}I_{i-2,2k+2}^+(h)) \\
	= &\sum\limits_{i+2k=0,i\geq 1}^{n-1}dA_{i+1,2k}^+I_{i,2k}^+(h)-\sum\limits_{i+2k=0,k\geq 1}{n}\frac{i+1}{2k+1}\frac{1}{c}A_{i+2,2k-2}^+I_{i,2k}^+(h) \\
	= & \sum\limits_{i+2k=0}^{n}\bar{B}_{i,2k}^+I_{1,2k}^+(h),
\end{aligned}
\end{equation*}
where
\begin{equation*}
\begin{aligned}
	\bar{B}_{i,2k}^+=L(A_{i,2k}^+|i+2k\leq n,i\geq 2,k\geq 0).
\end{aligned}	
\end{equation*}
Consequently, we get that
\begin{equation*}
	I_{n+1}^+(h)=\sum\limits_{k=0}^{[\frac{n}{2}]}A_{0,2k}^+I_{0,2k}^+(h)+\sum\limits_{k=0}^{[\frac{n-1}{2}]}A_{1,2k}^+(dI_{0,2k}^+(h)-\frac{2}{2k+3}\frac{1}{c}(2h)^{k+\frac{3}{2}})+\sum\limits_{i+2k=0}^n\tilde{A}_{i,2k}^+I_{i,2k}^+(h).
\end{equation*}
Hence, the expression of $I_n^+(h)$ in (\ref{eq_M_n(h)}) is proved by mathematical induction.

In a similar way, differentialing both sides of $H^-(x,y)=h$, we obtain that
\begin{equation*}
	ydy-ax^3dx+3abx^2dx-3ab^2xdx+ab^3dx=0.
\end{equation*}
Then, we get the following equation by multiplying both sides of it by $x^{i-3}y^{2k+1}$ and integrating it along $L_h^-$,
\begin{equation}
\begin{aligned}
	a\int\limits_{L_h^-}x^iy^{2k+1}dx-3ab\int\limits_{L_h^-}x^{i-1}y^{2k+1}dx+3ab^2\int_{L_h^-}x^{i-2}y^{2k+1}dx \\
	-ab^3\int\limits_{L_h^-}x^{i-3}y^{2k+1}dx-\int\limits_{L_h^-}x^{i-3}y^{2k+2}dy=0.
\end{aligned}
\end{equation}
As a result, we have that
\begin{equation}
\begin{aligned}
	& I_{i,2k}^-(h)=3bI_{i-1,2k}^-(h)-3b^2I_{i-2,2k}^-(h)+b^3I_{i-3,2k}^-(h)-\frac{i-3}{2k+3}\frac{1}{a}I_{i-4,2k+2}^-(h),~i\geq 4, \\
	& I_{3,2k}^-(h)-3b^2I_{1,2k}^-(h)+b^3I_{0,2k}^-(h)-\frac{2}{2k+3}\frac{1}{a}(2h)^{k+\frac{3}{2}}.
\end{aligned}
\end{equation}
Finally, we can also obtain the expression of $I_n^-(h)$ in~(\ref{eq_M_n(h)}). So it ends the proof.
\end{proof}

Considering with Lemmas~\ref{lemma1}-\ref{lemma3}, we have the following lemma.

\begin{lemma}\label{lemma4}
	Taking into account with~(\ref{eq_H}) and~(\ref{eq_L_h}), for $h\in(0,\frac{1}{4}ab^4)$, one has that
	
		\item $I^+(h)$ in~(\ref{M_h}) has the following form that
		\begin{equation}\label{M_hB}
			I^+(h)=\sum\limits_{k=0}^{[\frac{n-1}{2}]}D_{0,2k}^+(2h-cd^2)^kC_{i0}^+(h)+\sum\limits_{k=0}^{[\frac{n}{2}]}D_{1,2k}^+(2h)^{k+\frac{1}{2}},~n\geq 1,
		\end{equation}
		where
		$C_{i0}^+(h)$ appeared in~(\ref{I_i}), and $D_{00}^+,D_{02}^+,\cdots,B_{0,2[\frac{n-1}{2}]}^+,D_{10}^+,D_{12}^+,\cdots,D_{1,2[\frac{n}{2}]}$ are independent with each other.
		\item $I^-(h)$ in~(\ref{M_h}) has the following expression that
		      \begin{equation*}
			  \begin{aligned}
			  	I^-(h)= & \sum\limits_{k=0}^{[\frac{n-1}{2}]}D_{0,2k}^-(2h-\frac{1}{2}ab^4)^kC_{i0}^-(h)+\sum\limits_{k=0}^{[\frac{n-2}{2}]}D_{1,2k}^-(2h-\frac{1}{2}ab^4)^kC_{i1}^-(h) \\
			  	& +\sum\limits_{k=0}^{[\frac{n-3}{2}]}D_{2,2k}^-(2h-\frac{1}{2}ab^4)^kC_{i2}^-(h)+\sum\limits_{k=0}^{[\frac{n}{2}]}D_{3,2k}^-(2h)^{k+\frac{1}{2}},~n\geq 3,
			  \end{aligned}
			  \end{equation*}
			  where $C_{i0}^-(h)$ are shown in~(\ref{I_i}) and $D_{00}^-,D_{02}^-,\cdots,D_{0,2[\frac{n-1}{2}]}^-$, $D_{10}^-,D_{12}^-, \cdots,D_{1,2[\frac{n-2}{2}]}^-$, $D_{20}^-,\\D_{22}^-,\cdots, D_{2,2[\frac{n-3}{2}]}^-$, $D_{30}^-,D_{32}^-,\cdots,D_{3,2[\frac{n}{2}]}^-$
			  can be seen as free parameters. Here
	\begin{equation*}
	\begin{aligned}
		& D_{00}^+=2c^{-\frac{1}{2}}A_{00}^++L(A_{i,2k}^+|i+2k\leq n-1,i\geq 1,k\geq 0), \\
		& D_{0,2k}^+=2c^{-\frac{1}{2}}\prod\limits_{j=1}^k\beta_0(j)A_{0,2k}^++L(A_{i,2k}^+|i+2k\leq n-1,i\geq 1,k\geq 0), \quad k=1,2,\cdots,[\frac{n-1}{2}],\\
		& D_{10}^+=\frac{2a_{0,2k}^+}{2k+1}+L(A_{i,2k}^+|i+2k\leq n-1,i\geq 1,k\geq 0),\quad k=0,2,\cdots,[\frac{n}{2}].
	\end{aligned}
	\end{equation*}
	
\end{lemma}

\begin{proof}
	For $n\geq 1$, from Lemma~\ref{lemma2} and~\ref{lemma3}, it follows that
	\begin{equation*}
	\begin{aligned}
		I^+(h)= & I_n^+(h)+\sum\limits_{k=0}^{[\frac{n}{2}]}\frac{2a_{0,2k}^+}{2k+1}(2h)^{k+\frac{1}{2}} \\
		= & \sum\limits_{k=0}^{[\frac{n-1}{2}]}B_{0,2k}^+I_{0,2k}^+(h)+\sum\limits_{k=0}^{[\frac{n-2}{2}]}B_{1,2k}^+(2h)^{k+\frac{3}{2}}+\sum\limits_{k=0}^{[\frac{n}{2}]}\frac{2a_{0,2k}^+}{2k+1}(2h)^{k+\frac{1}{2}},
		\end{aligned}
	\end{equation*}
	which denotes as
	\begin{equation}\label{M_h_hat}
		I^+(h)=\sum\limits_{k=0}^{[\frac{n-1}{2}]}B_{0,2k}^+I_{0,2k}^+(h)+\sum\limits_{k=0}^{[\frac{n}{2}]}\hat{B}_{1,2k}^+(2h)^{k+\frac{1}{2}},~n\geq 1,
	\end{equation}
	where
	\begin{equation*}
		\hat{B}_{1,2k}^+=\frac{2a_{0,2k}^+}{2k+1}+L(B_{1,2k}^+|k=0,1,\cdots,[\frac{n-3}{2}]).
	\end{equation*}
	Noticing that $L_h^+$ is divided into two parts by $x$-axis which can be written as
	\begin{equation*}
	\begin{aligned}
		& y=(2h+c(x-d)^2-cd^2)^\frac{1}{2},~x\in(0,d-\sqrt{\frac{cd^2-2h}{c}}), \\
		& y=-(2h+c(x-d)^2-cd^2)^\frac{1}{2},~x\in(0,d-\sqrt{\frac{cd^2-2h}{c}}),
		\end{aligned}
	\end{equation*}
	we have 
	\begin{equation*}
		I_{i,2k}^+(h)=\int\limits_{L_h^+}x^iy^{2k+1}dx=2\int\limits_{0}^{d-\sqrt{\frac{cd^2-2h}{c}}}x^i(2h+c(x-d)^2-cd^2)^{k+\frac{1}{2}}dx.
	\end{equation*}
Let $t=-c^\frac{1}{2}(x-d)$ and then one has that
	\begin{equation}\label{I_1,2k+}
		I_{i,2k}^+(h)=2c^{-\frac{1}{2}}\int\limits_{\sqrt{cd^2-2h}}^{c^\frac{1}{2}d}(d-tc^{-\frac{1}{2}})^i(2h-cd^2+t^2)^{k+\frac{1}{2}}dt.
	\end{equation}
Based on ~(\ref{M_h_hat}) and (\ref{I_1,2k+}), it follows that
	\begin{equation*}
	\begin{aligned}
		I^+(h) & =\sum\limits_{k=0}^{[\frac{n}{2}]}\hat{B}_{1,2k}^+(2h)^{k+\frac{1}{2}}+\sum\limits_{k=0}^{[\frac{n-1}{2}]}2c^{-\frac{1}{2}}B_{0,2k}^+\int\limits_{\sqrt{cd^2-2h}}^{c^\frac{1}{2}d}(2h-cd^2+t^2)^{k+\frac{1}{2}}dt \\
		& =\sum\limits_{k=0}^{[\frac{n}{2}]}\hat{B}_{1,2k}^++\sum\limits_{k=0}^{[\frac{n-1}{2}]}\hat{B}_{0,2k}^+\int\limits_{\sqrt{cd^2-2h}}^{c^\frac{1}{2}d}(2h-cd^2+t^2)^{k+\frac{1}{2}}dt,
	\end{aligned}
	\end{equation*}
	where
	\begin{equation*}
	\begin{aligned}
		\hat{B}_{0,2k}^+ & =2c^{-\frac{1}{2}}A_{0,2k}^++L(A_{i,2k}^+|i+2k\leq n-1,i\geq 1,k\geq 0), \\
		\hat{B}_{1,2k}^+ & =\frac{2a_{0,2k}^+}{2k+1}+L(A_{i,2k}^+|i+2k\leq n-1,i\geq 1,k\geq 0).
	\end{aligned}
	\end{equation*}
Hence, from Lemma~\ref{lemma1}, for $n\geq 5$, we derive that 
	\begin{equation*}
	\begin{aligned}
		I^+(h) & =\sum\limits_{k=0}^{[\frac{n}{2}]}\hat{B}_{1,2k}^+(2h)^{k+\frac{1}{2}}+\hat{B}_{00}^+C_{i0}^+(h)+\hat{B}_{02}^+(\rho_0(1)(2h)^\frac{3}{2}+\bar{\rho}_0(1)(2h-cd^2)C_{i0}^+(h))\\ 
		&+  \sum\limits_{k=2}^{[\frac{n-1}{2}]}\hat{B}_{0,2k}^+(\alpha_0(k)(2h)^{k+\frac{1}{2}}+\sqrt{2h}\sum\limits_{j=1}^{k-1}\alpha_0(k-j)\prod\limits_{l=0}^{j-1}\beta_0(k-l)(2h-cd^2)^j(2h)^{k-j}\\&+\prod\limits_{j=1}^k\beta_0(j)(2h-cd^2)^kC_{i0}^+(h)).
	\end{aligned}
	\end{equation*}
which means (\ref{M_hB}) holds.
Denote matrix
	\begin{equation*}
	\begin{aligned}
		\mathcal{A} & =\frac{\partial (D_{00}^+,D_{02}^+,\cdots,D_{0,2[\frac{n-1}{2}]}^+,D_{10}^+,D_{12}^+,\cdots,D_{1,2[\frac{n}{2}])}^+)}{\partial (A_{00}^+,A_{02}^+,\cdots,A_{0,2[\frac{n-1}{2}]},a_{00}^+,a_{02}^+,\cdots,a_{0,2[\frac{n}{2}]})} =
		\begin{pmatrix}
			\mathcal{A}_{11} & 0 \\
			\mathcal{A}_{21} & \mathcal{A}_{22} \\
		\end{pmatrix},
	\end{aligned}
	\end{equation*}
	where $\mathcal{A}_{11},\mathcal{A}_{21},\mathcal{A}_{22}$ are $[\frac{n+1}{2}]\times[\frac{n+1}{2}],[\frac{n+2}{2}]\times[\frac{n+1}{2}],[\frac{n+2}{2}]\times[\frac{n+2}{2}]$ matrices respectively and
	\begin{equation*}
	\begin{aligned}
		\mathcal{A}_{11} & =2c^{-\frac{1}{2}}\mathbf{diag}(1,\beta_0(1),\prod_{j=1}^{2}\beta_0(j),\cdots,\prod\limits_{j=1}^{[\frac{n-1}{2}]}\beta_0(j)), \\
		\mathcal{A}_{22} & =2\mathbf{diag}(1,\frac{1}{3},\frac{1}{5},\cdots,\frac{1}{2[\frac{n}{2}]+1}).
	\end{aligned}
	\end{equation*}
Therefore, we have
	\begin{equation*}
	\begin{aligned}	
	|\mathcal{A}| & =|\mathcal{A}_{11}||\mathcal{A}_{22}|  =2^{[\frac{n+1}{2}]+[\frac{n+2}{2}]}c^{-\frac{1}{2}[\frac{n+1}{2}]}\prod\limits_{k=1}^{[\frac{n}{2}]}\frac{1}{2k+1}\prod\limits_{k=1}^{[\frac{n-1}{2}]}(\prod\limits_{j=1}^{k}\beta_0(j))  \neq 0.
	\end{aligned}
	\end{equation*}
We can easily know the coefficients in~(\ref{M_hB}) are free parameters.  Now, we prove (ii).
	
	Analogously, for $n\geq 7$, one has that
	\begin{equation*}
	\begin{aligned}
		I^-(h)= & \sum\limits_{k=0}^{[\frac{n-1}{2}]}B_{0,2k}^-I_{0,2k}^-(h)+\sum\limits_{k=0}^{[\frac{n-1}{2}]}B_{1,2k}^-I_{1,2k}^-(h) \\
		& +\sum\limits_{k=0}^{[\frac{n-3}{2}]}B_{2,2k}^-I_{2,2k}^-(h)+\sum\limits_{k=0}^{[\frac{n-4}{2}]}B_{3,2k}^-(2h)^{k+\frac{3}{2}}+\sum\limits_{k=0}^{[\frac{n}{2}]}\frac{2a_{0,2k}^-}{2k+1}(2h)^{k+\frac{1}{2}}.
	\end{aligned}
	\end{equation*}
	Then, we have
	\begin{equation}\label{M^-hhat}
	\begin{aligned}
	I^-(h)= & \sum\limits_{k=0}^{[\frac{n-1}{2}]}B_{0,2k}^-I_{0,2k}^-(h)+\sum\limits_{k=0}^{[\frac{n-1}{2}]}B_{1,2k}^-I_{1,2k}^-(h)  +\sum\limits_{k=0}^{[\frac{n-3}{2}]}B_{2,2k}^-I_{2,2k}^-(h)+\sum\limits_{k=0}^{[\frac{n}{2}]}\hat{B}_{3,2k}^-(2h)^{k+\frac{1}{2}}
	\end{aligned}
	\end{equation}
	where
	\begin{equation*}
		\hat{B}_{3,2k}^-=\frac{2a_{0,2k}^-}{2k+1}+L(\hat{B}_{3,2k}^-|k=0,1,\cdots,[\frac{n-4}{2}]).
	\end{equation*}
	Then, letting $t=-(\dfrac{a}{2})^\frac{1}{4}(x-b)$, one has that
	\begin{equation*}
	\begin{aligned}
		I_{i,2k}^-(h) & =\int\limits_{L_h^-}x^iy^{2k+1}dx=2\int\limits_{{b+(\frac{ab^4-4h}{a})}^\frac{1}{4}}^0x^i(2h+\frac{1}{2}a(x-b)^4-\frac{1}{2}ab^4)^{k+\frac{1}{2}}dx \\ =&2(\frac{2}{a})^\frac{1}{4}\int\limits_{{(\frac{a}{2})}^\frac{1}{4}b}^{-(\frac{ab^4}{2}-2h)^\frac{1}{4}}(b-(\frac{2}{a})^\frac{1}{4}t)^i(2h+t^4-\frac{1}{2}ab^4)^{k+\frac{1}{2}}dt.
	\end{aligned}
	\end{equation*}\vspace*{-0.2cm}Hence,(\ref{M^-hhat}) can be rewritten as
	\begin{equation}
	\begin{aligned}
		I^-(h)= & \sum\limits_{k=0}^{[\frac{n-1}{2}]}2(\frac{2}{a})^\frac{1}{4}B_{0,2k}^-\int\limits_{(\frac{a}{2})^\frac{1}{4}b}^{-(\frac{ab^4}{2}-2h)^\frac{1}{4}}(2h+t^4-\frac{1}{2}ab^4)^{k+\frac{1}{2}}dt \\
		& +\sum\limits_{k=0}^{[\frac{n-2}{2}]}2(\frac{2}{a})^\frac{1}{4}B_{1,2k}^-\int\limits_{(\frac{a}{2})^\frac{1}{4}b}^{-(\frac{ab^4}{2}-2h)^\frac{1}{4}}(b-(\frac{2}{a})^\frac{1}{4}t)(2h+t^4-\frac{1}{2}ab^4)^{k+\frac{1}{2}}dt \\
		& +\sum\limits_{k=0}^{[\frac{n-3}{2}]}2(\frac{2}{a})^\frac{1}{4}B_{2,2k}^-\int\limits_{(\frac{a}{2})^\frac{1}{4}b}^{-(\frac{ab^4}{2}-2h)^\frac{1}{4}}(b-(\frac{2}{a})^\frac{1}{4}t)^2(2h+t^4-\frac{1}{2}ab^4)^{k+\frac{1}{2}}dt  +\sum\limits_{k=0}^{[\frac{n}{2}]}\hat{B}_{3,2k}^-(2h)^{k+\frac{1}{2}} \\
		= & \sum\limits_{k=0}^{[\frac{n-1}{2}]}\hat{B}_{0,2k}^-\int\limits_{(\frac{a}{2})^\frac{1}{4}b}^{-(\frac{ab^4}{2}-2h)^\frac{1}{4}}(2h+t^4-\frac{1}{2}ab^4)^{k+\frac{1}{2}}dt \\
		& +\sum\limits_{k=0}^{[\frac{n-2}{2}]}\hat{B}_{1,2k}^-\int\limits_{(\frac{a}{2})^\frac{1}{4}b}^{-(\frac{ab^4}{2}-2h)^\frac{1}{4}}(b-(\frac{2}{a})^\frac{1}{4}t)(2h+t^4-\frac{1}{2}ab^4)^{k+\frac{1}{2}}dt \\
		& +\sum\limits_{k=0}^{[\frac{n-3}{2}]}\hat{B}_{2,2k}^-\int\limits_{(\frac{a}{2})^\frac{1}{4}b}^{-(\frac{ab^4}{2}-2h)^\frac{1}{4}}(b-(\frac{2}{a})^\frac{1}{4}t)^2(2h+t^4-\frac{1}{2}ab^4)^{k+\frac{1}{2}}dt \\
		& +\sum\limits_{k=0}^{[\frac{n}{2}]}\hat{B}_{3,2k}^-(2h)^{k+\frac{1}{2}},
	\end{aligned}
	\end{equation}
	where
	\begin{equation*}
	\begin{aligned}
		 & \hat{B}_{0,2k}^-=2(\frac{2}{a})^\frac{1}{4}A_{0,2k}^-+L(A_{i,2k}^-|i+2k\leq n-1,i\geq 1,k\geq 0), \\
		 & \hat{B}_{1,2k}^-=-2(\frac{2}{a})^\frac{1}{2}A_{1,2k}^-+L(A_{i,2k}^-|i+2k\leq n-1,i\geq 2,k\geq 0), \\	
		 & \hat{B}_{2,2k}^-=2(\frac{2}{a})^\frac{3}{4}A_{2,2k}^-+L(A_{i,2k}^-|i+2k\leq n-1,i\geq 3,k\geq 0), \\	
		 & \hat{B}_{3,2k}^-=\frac{2a_{0,2k}^-}{2k+1}+L(A_{i,2k}^-|i+2k\leq n-1,i\geq 3,k\geq 0).
	\end{aligned}
	\end{equation*}
From Lemma~\ref{lemma1}, we derive the conclusion in (ii). This finishes the proof.
\end{proof}

Next, from Lemma~\ref{lemma4}, we have the following proposition.
\begin{proposition}\label{prop1}
	Assume that~(\ref{eq_pw}) and~(\ref{eq_L_h}) hold. Then, for $n\geq 1$, $I(h)$ has the following expression when $h\in(0,\dfrac{1}{2}\rho)$,
	\begin{equation*}
	\begin{aligned}
		I(h)= & \sum\limits_{k=0}^{[\frac{n}{2}]}A_{0,2k}(2h)^{k+\frac{1}{2}}+\sum\limits_{k=0}^{[\frac{n-1}{2}]}A_{1,2k}(2h-\rho)^kC_{i0}^+(h) \\
		& +\sum\limits_{k=0}^{[\frac{n-1}{2}]}A_{2,2k}(2h-\rho)^kC_{i0}^-(h)+\sum\limits_{k=0}^{[\frac{n-2}{2}]}A_{3,2k}(2h-\rho)^kC_{i1}^-(h) \\
		& +\sum\limits_{k=0}^{[\frac{n-3}{2}]}A_{4,2k}(2h-\rho)^kC_{i2}^-(h),
	\end{aligned}
	\end{equation*}
	where $\rho =\dfrac{1}{4}ab^4$, $A_{0,2k}(k=0,\cdots,[\frac{n}{2}])$, $A_{l,2k}(l=1,2,k=0,\cdots,[\frac{n-1}{2}])$, $A_{3,2k}(k=0,\cdots,[\frac{n-2}{2}])$, $A_{4,2k}(k=0,\cdots,[\frac{n-3}{2}])$ can be seen as free parameters.
\end{proposition}
	
This then implies the following corollary:
\begin{corollary}
	Suppose that~(\ref{eq_pw}) and~(\ref{eq_L_h}) hold. Then, $I(h)$ has the following expression for $h\in (0,\dfrac{1}{\rho})$,
	\begin{equation*}
		I(h)=
		\begin{cases}
			A_{00}(2h)^\frac{1}{2}+A_{10}C_{i0}^+(h)+A_{20}C_{i0}^-(h),\quad n=1, \\
			A_{00}(2h)^\frac{1}{2}+A_{02}(2h)^\frac{3}{2}+A_{10}C_{i0}^+(h)+A_{20}C_{i0}^-(h)+A_{30}C_{i1}^-(h),\quad n=2,
		\end{cases}
	\end{equation*}
	where the coefficients are independent with each other.
\end{corollary}

\subsection{The Asymptotic Expansion of $I(h)$ near $L$}

Firstly, we denote
\begin{equation}\label{u}
	u=2h-\rho,
\end{equation}
and introduce the following lemma.
\begin{lemma}\label{lemma7}
	For $0<-u\ll 1$, $C_{i0}^+(h)$ can be expressed as
	\begin{equation}\label{I_0h2}
		C_{i0}^+(h)=\eta_0u\ln|u|+\sum\limits_{k\geq 0}\eta_{0k}u^k,\eta_0<0.
	\end{equation}
where $\eta_0$, $\eta_{0k},k\geq 0$ are all real constants.
\end{lemma}

\begin{proof}
We have the expression of $I_i^+(h)$ with~(\ref{u}) that
	\begin{equation}
		C_{i0}^+(h)=\int\limits_{(cd^2-2h)^{\frac{1}{2}}}^{c^{\frac{1}{2}}d}t^i(2h+t^2-cd^2)^{\frac{1}{2}}dt=\int\limits_{|u|^\frac{1}{2}}^{c^\frac{1}{2}d}t^i(u+t^2)^\frac{1}{2}dt=\int\limits_{|u|^\frac{1}{2}}^{c^\frac{1}{2}d}t^{i+1}(ut^{-2}+1)^\frac{1}{2}dt.
	\end{equation}
Based on $t\in [(\frac{a}{2})^\frac{1}{4}b,-|u|^\frac{1}{4}]$, we have $ut^{-4}\in [-1,u(\frac{a}{2})^{-1}b^{-4}]\subset [-1,0)$.
	Hence, $C_{i0}^+(h)$ can be written as
	\begin{equation*}
	\begin{aligned}
		C_{i0}^+(h) & =\int\limits_{|u|^\frac{1}{2}}^{c^\frac{1}{2}d}t^{i+1}(ut^{-2}+1)^\frac{1}{2}dt  = \int\limits_{|u|^\frac{1}{2}}^{c^\frac{1}{2}d}t^{i+1}\sum\limits_{k\geq 0}d_ku^kt^{-2k} =\sum\limits_{k\geq 0}d_ku^k\int\limits_{|u|^\frac{1}{2}}^{c^\frac{1}{2}d}t^{i-2k+1}dt \\ &=\sum\limits_{k\geq 0}d_ku^k\frac{1}{i-2k+2}\left.t^{i-2k+2}\right|_{|u|^\frac{1}{2}}^{c^\frac{1}{2}d},
	\end{aligned}
	\end{equation*}
	where $d_k$, $k\geq 0$ satisfy $(1+x)^\frac{1}{2}=\sum\limits_{k\geq 0}d_kx^k$, $x\in [-1,1]$.
	By calculating, it is not hard to get 
	\begin{equation}
		d_0=1,~d_1=\frac{1}{2},~d_k=\frac{(-1)^{k-1}(2k-3)!!}{(2k)!!},~k\geq 2.
	\end{equation}
	Then, it follows that
	\begin{equation*}
	\begin{aligned}
		C_{i0}^+(h) & =\sum\limits_{k\geq 0,k\neq 1}d_ku^k\frac{1}{2-2k}\left.t^{2-2k}\right|_{|u|^\frac{1}{2}}^{c^\frac{1}{2}d}+d_1u\left.\ln t\right|_{|u|^\frac{1}{2}}^{c^\frac{1}{2}d} \\
		& =\sum\limits_{k\geq 0,k\neq 1}\frac{d_k}{2-2k}c^{1-k}d^{2-2k}u^k+\sum\limits_{k\geq 0,k\neq 1}\frac{(-1)^kd_k}{2k-2}|u|+d_1u\ln (c^\frac{1}{2}d)-\frac{1}{2}d_1u\ln |u|.
	\end{aligned}
	\end{equation*}
Immediately, we obtain the expression in~(\ref{I_0h2}) by denoting
\begin{equation}
\begin{aligned}
	& \eta_{0k}=\frac{d_k}{2-2k}c^{1-k}d^{2-2k},~k\geq 0,~k\neq 1 \\
	& \eta_{0}=-\frac{1}{2}d_1,
\end{aligned}
\end{equation}
where
\begin{equation*}
	\delta_0<0.
\end{equation*}
This ends the proof.
\end{proof}

\begin{lemma}\label{lemma8}
	For $0<-u\ll 1$, we have
	\begin{equation}
	\begin{aligned}
		& C_{i0}^-(h)=\sigma_0|u|^\frac{3}{4}+\sum\limits_{k\geq 0}\sigma_{0k}u^k,~C_{i1}^-(h)=\sigma_1u\ln|u|+\sum\limits_{k\geq 0}\sigma_{1k}u^k,  ~C_{i2}^-(h)=\sigma_2|u|^\frac{5}{4}+\sum\limits_{k\geq 0}\sigma_{2k}u^k,
	\end{aligned}
	\end{equation}
	where
	\begin{equation*}
		\sigma_0< 0,\quad \sigma_1 <0,\quad \sigma_2 >0.
	\end{equation*}
\end{lemma}

\begin{proof}
	Similar to the proof of Lemma~\ref{lemma7}, one has that
	\begin{equation*}
	C_{i0}^-(h)=\int\limits_{(\frac{a}{2})^\frac{1}{4}b}^{-(\frac{ab^4}{2}-2h)^\frac{1}{4}}t^i(2h+t^4-\frac{1}{2}ab^4)^\frac{1}{2}dt=\int\limits_{(\frac{a}{2})^\frac{1}{4}b}^{-|u|^\frac{1}{4}}t^i(u+t^4)^\frac{1}{2}dt=\int\limits_{(\frac{a}{2})^\frac{1}{4}b}^{-|u|^\frac{1}{4}}t^{i+2}(ut^{-4}+1)^\frac{1}{2}dt.
	\end{equation*}
	Since $t\in[(\frac{a}{2})^\frac{1}{4}b,-|u|^\frac{1}{4}]$, we have $ut^{-4}\in [-1,u(\frac{a}{2})^{-1}b^{-4}]\subset[-1,0)$.
	Similarily, we get
	\begin{equation*}
	\begin{aligned}
		C_{i0}(h) & =\int\limits_{(\frac{a}{2})^\frac{1}{4}b}^{-|u|^\frac{1}{4}}t^{i+2}\sum\limits_{k\geq 0}d_ku^kt^{-4k}dt \\
		& =\sum\limits_{k\geq 0}d_ku^k\int\limits_{(\frac{a}{2})^\frac{1}{4}b}^{-|u|^\frac{1}{4}}t^{i+2-4k}dt \\
		& =
		\begin{cases}
			 \sum\limits_{k\geq 0}d_ku^k\frac{1}{i+3-4k}\left.t^{i+3-4k}\right|_{(\frac{a}{2})^\frac{1}{4}b}^{-|u|^\frac{1}{4}},\quad i=0,2, \\
		     \sum\limits_{k\geq 0,k\neq 1}d_ku^k\frac{1}{4-4k}\left.t^{4-4k}\right|_{(\frac{a}{2})^\frac{1}{4}b}^{-|u|^\frac{1}{4}}+d_1u\left.\ln -t\right|_{(\frac{a}{2})^\frac{1}{4}b}^{-|u|^\frac{1}{4}},\quad i=1,
		\end{cases}
		\\
		& =
		\begin{cases}
			\sum\limits_{k\geq 0}\frac{(-1)^{k+1}d_k}{i+3-4k}|u|^\frac{i+3}{4}-\sum\limits_{k\geq 0}\frac{d_k}{i+3-4k}(\frac{a}{2})^\frac{i+3-4k}{4}b^{i+3-4k}u^k,\quad i=0,2 \\
			\sum\limits_{k\geq 0,k\neq 1}\frac{(-1)^kd_k}{4-4k}|u|-\sum\limits_{k\geq 0,k\neq 1}\frac{\xi_k}{4-4k}(\frac{a}{2})^{1-k}b^{4-4k}u^k+\frac{1}{4}d_1u\ln |u|-d_1u\ln ((\frac{a}{2})^\frac{1}{4}b),\quad i=1,
		\end{cases}
	\end{aligned}
	\end{equation*}
	where
	\begin{equation*}
	\begin{aligned}
		& \sigma_{ik}=\frac{d_k}{4k-i-3}(\frac{a}{2})^\frac{i+3-4k}{4}b^{i+3-4k},\quad i=0,2,\quad k\geq 0, \\
		& \sigma_{1k}=\frac{d_k}{4k-4}(\frac{a}{2})^{1-k}b^{4-4k},\quad k\geq 0,\quad k\neq 1, \\
		& \sigma_{11}=\sum\limits_{k\geq 0,k\neq 1}\frac{(-1)^{k+1}d_k}{4-4k}-d_1\ln (\frac{a}{2})^\frac{1}{4}b, \\
		& \sigma_i=\sum\limits_{k\geq 0}\frac{(-1)^{k+1}d_k}{i+3-4k},\quad i=0,2, \\
		& \sigma_1=\frac{1}{4}d_1,
	\end{aligned}
	\end{equation*}
	with
	\begin{equation*}
		\sigma_0< 0,\quad\sigma_1<0,\quad\sigma_2>0.
	\end{equation*}
\end{proof}

Then, it follows that

\begin{proposition}\label{prop2}
	Assume that~(\ref{eq_H}) and~(\ref{eq_L_h}) hold, for $n\geq 1$, $I(h)$ near  $L$ have the following expression:
	\begin{equation*}
		I(h)=c_{10}+c_{30}|u|^\frac{3}{4}+\sum\limits_{i=1}^{[\frac{n+1}{2}]}(c_{0i}\ln |u|+c_{1i}+c_{2i}|u|^\frac{1}{4}+c_{3i}|u|^\frac{3}{4})|u|^i+\sum\limits_{i\geq [\frac{n+1}{2}]+1}c_{1i}|u|^i,\quad 0<-u\ll 1,
	\end{equation*}
	the coefficients $c_{0i},i=1,2,\cdots,[\frac{n+1}{2}]$, $c_{1i},i=0,1,\cdots,[\frac{n+1}{2}]$, $c_{2i},i=1,2,\cdots,[\frac{n-1}{2}]$, $c_{3i},i=0,1,\cdots,[\frac{n-1}{2}]$ can be taken as free parameters.
	Here
	\begin{equation*}
	\begin{aligned}
		c_{i}= & L(c_{0i}|i=1,2,\cdots,[\frac{n+1}{2}])+L(c_{1i}|i=0,1,\cdots,[\frac{n+1}{2}]) \\
		& +L(c_{2i}|i=1,2,\cdots,[\frac{n-1}{2}])+L(c_{3i}|i=1,2,\cdots,[\frac{n-1}{2}]),i\geq [\frac{n+1}{2}]+1.
	\end{aligned}
	\end{equation*}	
\end{proposition}

\begin{proof}
	From Proposition~\ref{prop1}, Lemma \ref{lemma7} and Lemma 3.6 in \cite{Xiong2021}, for $n\geq 1$, $0<-u\ll 1$, we have
	\begin{equation*}
	\begin{aligned}
		I(h)= & \sum\limits_{i\geq 0}c_iu^i+\sum\limits_{k=0}^{[\frac{n-1}{2}]}A_{1,2k}u^k(\eta_0u\ln|u|+\sum\limits_{k\geq 0}\eta_{0k}u^k) \\
		& +\sum\limits_{k=0}^{[\frac{n-2}{2}]}A_{2,2k}u^k(\sigma_0|u|^\frac{3}{4}+\sum\limits_{k\geq 0}\sigma_{0k}u^k) \\
		& + \sum\limits_{k=0}^{[\frac{n-3}{2}]}A_{4,2k}u^k(\sigma_2|u|^\frac{5}{4}+\sum\limits_{k\geq 0}\sigma_{2k}u^k).
	\end{aligned}
	\end{equation*}
	Then, one has
	\begin{equation*}
	\begin{aligned}
		I(h)= & \sum\limits_{i=1}^{[\frac{n+1}{2}]}c_{0i}|u|^i\ln |u|+\sum\limits_{i\geq 0}c_{1i}|u|^i  +\sum\limits_{i=1}^{[\frac{n-1}{2}]}c_{2i}|u|^{\frac{1}{4}+i}+\sum\limits_{i=0}^{[\frac{n-1}{2}]}c_{3i}|u|^{\frac{3}{4}+i},
	\end{aligned}
	\end{equation*}
	which is the expression in Proposition\ref{prop2}. The coefficients are
	\begin{equation*}
	\begin{aligned}
		c_{0i}= & (-1)^i\eta_0A_{1,2(i-1)}+(-1)^i\sigma_1A_{3,2(i-1)},i=1,2,\cdots,[\frac{n+1}{2}], \\
		c_{1i}= & c_i+L(A_{1,2k}|k=0,\cdots,[\frac{n-1}{2}])+L(A_{2,2k}|k=0,\cdots,[\frac{n-1}{2}]) \\
		& +L(A_{3,2k}|k=0,\cdots,[\frac{n-2}{2}])+L(A_{4,2k}|k=0,\cdots,[\frac{n-3}{2}]),i=0,1,\cdots, \\
		c_{2i}= & (-1)^{i-1}\sigma_2A_{4,2(i-1)},i=1,\cdots,[\frac{n-1}{2}], \\
		c_{3i}= & (-1)^i\sigma_0A_{2,2i},i=0,1,\cdots,[\frac{n-1}{2}].
	\end{aligned}
	\end{equation*}
	From Lemma~\ref{lemma7} and proposition~\ref{prop1}, We can define a matrix $\mathcal{B}$
	\begin{equation*}
	\begin{aligned}
	\mathcal{B} & =\frac{\partial (c_{01},\cdots,c_{0,[\frac{n+1}{2}]},c_{10},\cdots,c_{1,[\frac{n+1}{2}]},c_{21},\cdots,c_{2,[\frac{n-1}{2}]},c_{31},\cdots,c_{3,[\frac{n-1}{2}]})}{\partial (A_{10},\cdots,A_{1,2[\frac{n-1}{2}]},c_0,\cdots,c_{[\frac{n+1}{2}]},A_{40},\cdots,A_{4,2[\frac{n-3}{2}]},A_{20},\cdots,A_{2,2[\frac{n-1}{2}]})} \\
	& =
	\begin{pmatrix}
	\mathcal{B}_{11} & 0 & 0 & 0 \\
	\mathcal{B}_{21} & \mathcal{B}_{22} & \mathcal{B}_{23} & \mathcal{B}_{24} \\
	0 & 0 & \mathcal{B}_{33} & 0 \\
	0 & 0 & 0 & \mathcal{B}_{44} \\
	\end{pmatrix}
	\end{aligned}
	\end{equation*}
	where
	\begin{equation*}
	\begin{aligned}
		& \mathcal{B}_{11}=\eta_0\mathbf{diag}(-1,1,\cdots,(-1)^{[\frac{n+1}{2}]}), \\
		& \mathcal{B}_{22}=\mathbf{diag}(1,1,\cdots,1), \\
		& \mathcal{B}_{33}=\sigma_2\mathbf{diag}(1,-1,\cdots,(-1)^{[\frac{n-3}{2}]}), \\
		& \mathcal{B}_{44}=\sigma_0\mathbf{diag}(1,-1,\cdots,(-1)^{[\frac{n-1}{2}]}).
	\end{aligned}
	\end{equation*}
	Hence, we have
	\begin{equation*}
	\begin{aligned}
		|\mathcal{B}|= & \prod_{i=1}^{4}|\mathcal{B}_{ii}|=\eta_0^{[\frac{n+1}{2}]}\sigma_0^{[\frac{n-1}{2}]}\sigma_2^{[\frac{n-1}{2}]}\prod\limits_{i=1}^{[\frac{n+1}{2}]}(-1)^i\prod\limits_{i=0}^{[\frac{n-3}{2}]}(-1)^i\prod\limits_{i=0}^{[\frac{n-1}{2}]}(-1)^i 
		=  -\eta_0^{[\frac{n+1}{2}]}\sigma_0^{[\frac{n-1}{2}]}\sigma_2^{[\frac{n-1}{2}]}\prod\limits_{i=1}^{[\frac{n-1}{2}]}(-1)^i.
	\end{aligned}
	\end{equation*}
	By Lemma~\ref{lemma7} and~\ref{lemma8}, we derive
	\begin{equation*}
		|\mathcal{B}|\neq 0.
	\end{equation*}
Therefore,  the coefficients referred in Proposition~\ref{prop2} are independent of each other. 
\end{proof}

\section{Proof of Theorem~\ref{thm1}}\label{proof}
In this section, we will give the proof of Theorem~\ref{thm1}.
\begin{proof}
From Proposition~\ref{prop2}, it is proved that the coefficients $c_{0i},i=1,2,\cdots,[\frac{n+1}{2}]$, $c_{1i},i=0,1,\cdots,[\frac{n+1}{2}]$, $c_{2i},i=1,2,\cdots,[\frac{n-1}{2}]$, $c_{3i},i=0,1,\cdots,[\frac{n-1}{2}]$ are all free parameters.
Let 
\begin{equation*}
\begin{aligned}
	& 0<c_{3,i-1}\ll c_{0i}\ll c_{1i}\ll -c_{2i}\ll c_{3i}, i=1,2,\cdots,[\frac{n+1}{2}], \\
	& 0< c_{10}\ll c_{30},\quad 0<c_{3,[\frac{n+1}{2}]}\ll 1.
\end{aligned}
\end{equation*}
such that the sigh of $I(h)$ in Proposition~\ref{prop2} has been changed $4[\frac{n+1}{2}]+1,n\geq 1$ times as $-u>0$ sufficiently small.
It implies that $I(h)$ can have $4[\frac{n+1}{2}]+1,n\geq 1$ simple zeros near $L$.
When
\begin{equation*}
\begin{aligned}
	c_{0i}=0,i=1,2,\cdots,[\frac{n+1}{2}],\\
	c_{1i}=0,i=0,1,\cdots,[\frac{n+1}{2}],\\
	c_{2i}=0,i=1,2,\cdots,[\frac{n-1}{2}],\\
	c_{3i}=0,i=0,1,\cdots,[\frac{n-1}{2}],
\end{aligned}
\end{equation*}
by Proposition~\ref{prop2}, we have
\begin{equation*}
	I(h)\equiv 0,h\in(0,\frac{1}{4}ab^4).
\end{equation*}
Therefore, it implies that the maximal number of simple zeros of $I(h)$ is $4[\dfrac{n+1}{2}]+1,n\geq 1$ as $-u>0$, i.e. $h<\dfrac{1}{4}ab^4$, is sufficiently small. 
\end{proof}

\section{Acknowledgments}
\noindent This work was supported by the National Natural Science Foundation of China (NNSFC) (No. 12172340), the Fundamental Research Funds for the Central Universities, China University of Geosciences (Wuhan) (Nos. CUGGC05 and CUGDCJJ202216), and the Young Top-notch Talent Cultivation Program of Hubei Province.

\bibliographystyle{ieeetr}
\bibliography{sample}
\end{document}